# Joint differential resolvents for pseudopolynomials

John Michael Nahay, Ph.D., Broadway Performance Systems, LLC Silver Spring, Maryland broadway-performance-systems.com

**Abstract** The existence of linear differential resolvents for $z^\alpha$ for any root $z$ of an ordinary polynomial with coefficients in a given ordinary differential field has been established, where $\alpha$ is an indeterminate constant with respect to the derivation of the given field.

In this paper we consider several alphas. We will call a finite sum of indeterminate powers of a variable $v$ a *pseudopolynomial* in $v$.

We will generalize the definition of a differential resolvent of a single polynomial for a single monomial $z^\alpha$ to the definition of a differential resolvent of several polynomials for a pseudopolynomial in the roots. We will also generalize the definition of a resolvent to have *non-consecutive* derivatives.

We will show that the author's powersum formula may be used to compute this more general differential resolvent.

## 1. Standard notation from differential algebra and Introduction

All rings will be commutative with unity 1.

Let $\mathbb{Z}$ denote the ring of integers.

Let $\mathbb{Z}^+$ denote the set of positive integers.

Let $\mathbb{Z}_0^+$ denote the set of nonnegative integers.

Let $\mathbb{Q}$ denote the field of rational numbers.

For any finite set $S$, let $|S|$ denote its size.

The symbol $\forall$ means "for all". The symbol $\ni$ means "such that".

The symbol $\equiv$ means "is defined as" or "is identically equal to".

For any positive integer, $m \in \mathbb{Z}^+$, define $[m] \equiv \{k \in \mathbb{Z}^+ \ni 1 \le k \le m\}$ and $[m]_0 \equiv \{k \in \mathbb{Z}_0^+ \ni 0 \le k \le m\}$.

For any ring $\mathbb{R}$, let $\mathbb{R}^\#$ denote the proper subset of nonzero elements of $\mathbb{R}$.

For any ring $\mathbb{R}$, let $\mathbb{R}[h_1,...,h_L]$ denote the ring generated by $\mathbb{R}$ and the variables $\{h_1,...,h_L\}$, that is, the set of all polynomials over $\mathbb{R}$ of $\{h_1,...,h_L\}$.

For any field $\mathbb{F}$, let $\mathbb{F}(h_1,...,h_L)$ denote the field generated by $\mathbb{F}$ and the variables $\{h_1,...,h_L\}$, that is, the set of all rational functions over $\mathbb{F}$ of the variables $\{h_1,...,h_L\}$.

A derivation $D$ on a ring $\mathbb{R}$ is a map from $\mathbb{R}$ to $\mathbb{R}$ which satisfies the Leibniz rule $D(u \cdot v) = u \cdot Dv + D(u) \cdot v$. If $r \in \mathbb{R}$, then $r \cdot D$ is also a derivation on $\mathbb{R}$. Hence, $0 \cdot D$, which maps $\mathbb{R}$ to 0, is a derivation. Henceforth, we will assume that our derivation $D$ is not the identically zero derivation, that is, $D$ is non-trivial.

A ring $\mathbb{R}$ equipped with a non-trivial derivation $D$ is called a differential ring.

For any differential ring $\mathbb{R}$, let $\mathbb{R}\{h_1,...,h_L\}$ denote the differential ring generated by $\mathbb{R}$ and the variables $\{h_1,...,h_L\}$, that is, the set of all polynomials over $\mathbb{R}$ of infinitely many derivatives of the variables $\{h_1,...,h_L\}$.

For any differential field $\mathbb{F}$, let $\mathbb{F}<h_1,...,h_L>$ denote the differential field generated by $\mathbb{F}$ and the variables $\{h_1,...,h_L\}$, that is, the set of all rational functions over $\mathbb{F}$ of infinitely many derivatives of the variables $\{h_1,...,h_L\}$.

Although nothing in Theorems 3.1 and 4.2 requires the differential fields to have zero characteristic, many other theorems on differential resolvents have been proven only for char($\mathbb{F}$)=0. Hence, unless we specify otherwise, assume zero characteristic for all differential fields. Example 1.5 will demonstrate a resolvent when char($\mathbb{F}$)=3.

Joseph Louis Lagrange published his famous formula for an inverse of a power series of a single variable in 1770 [5]. In Section 5.1 of [2] Egorychev generalized Lagrange's formula to a formula for multiple variables. Suppose $F_i(w, z_1,...,z_n) = 0$ are $n$ holomorphic functions of the $n$ variables $\{z_1,...,z_n\}$ and variables $w$. There is no loss of generality if we lump all the free variables together and call them by the single letter name $w$. Egorychev gave increasingly more explicit formulae for $\Phi(w, \varphi(w))$ where $\Phi$ is a given holomorphic function of its arguments, and $(z_1,...,z_n) = \varphi(w) = (\varphi_1(w),...,\varphi_n(w))$ are implicit functions determined by $F_i(w, z_1,...,z_n) = 0$.

When the $F_i$ are *polynomials* in the $n$ variables $\{z_1,...,z_n\}$, then there exists a linear differential resolvent – basically, a finite-order linear differential operator – for $\Phi$ whose terms depend upon $w$. We will coin a new term for this differential operator for $\Phi$ -- we will call it a *joint* differential resolvent *of* the polynomials $F_i$ *for* $\Phi$.

For any constant $\alpha$, transcendental or algebraic over $\mathbb{Q}<z_1,...,z_n>$, an $\alpha$-resolvent of a univariate polynomial $P(t)$ with $n$ roots $\{z_1,...,z_n\}$ is the special case of a differential resolvent of $P(t)$ for the pseudomonomial $\Phi = z^\alpha$.

To be consistent with notation on earlier documents on differential resolvents rather than Egorychev's notation, we will use the letter $y$ instead of $\Phi$ henceforth.

This article will generalize the definition of a differential $\alpha$-resolvent into *three* directions.

**Generalization 1 of 3.** First, we will generalize the definition of a differential resolvent of a single polynomial to the joint differential resolvent of several polynomials.

**Example 1.1.** Let $P_1(t)$ be a quadratic polynomial with roots $z_1 = e^x$ and $z_2 = \ln x$. In other words, $P_1(t) = (t - e^x) \cdot (t - \ln x) = t^2 - (e^x + \ln x) \cdot t + e^x \cdot \ln x$. Define a derivation $D$ such that $Dx \equiv 1$. Let $\mathbb{F}$ denote the smallest differential field generated by the coefficients of $P_1(t)$, in other words, $\mathbb{F} = \mathbb{Q}<e^x + \ln x, e^x \cdot \ln x>$. We call $\mathbb{F}$ the *differential coefficient field of* $P_1(t)$. By definition, a differential $\alpha$-resolvent of $P_1(t)$ is any nonzero finite-order linear differential operator $\mathfrak{R}$ whose terms lie in $\mathbb{F}<\alpha>$ such that $\mathfrak{R} e^{\alpha \cdot x} = 0$ and $\mathfrak{R}(\ln x)^\alpha = 0$.

Note that since $e^x \cdot \ln x \in \mathbb{F}$, then $D(e^x \cdot \ln x) \in \mathbb{F}$. So $e^x \cdot \ln x + e^x \cdot x^{-1} \in \mathbb{F}$. So $e^x \cdot x^{-1} \in \mathbb{F}$. Since $e^x + \ln x \in \mathbb{F}$, then $D(e^x + \ln x) \in \mathbb{F}$. So $e^x + x^{-1} \in \mathbb{F}$. The inclusions $e^x \cdot x^{-1} \in \mathbb{F}$ and $e^x + x^{-1} \in \mathbb{F}$ imply $x^{-2} \in \mathbb{F}$, so $x^2 \in \mathbb{F}$. So $D(x^2) \in \mathbb{F}$. So $2x \in \mathbb{F}$. So $x \in \mathbb{F}$. Note also that $\alpha$ constant implies $\mathbb{F}<\alpha> = \mathbb{F}(\alpha)$, and $\mathbb{R}\{\alpha\} = \mathbb{R}[\alpha]$ for any differential ring $\mathbb{R}$.

Since $\Re z_1^\alpha = 0$ and $\Re z_2^\alpha = 0$, then $\Re(a z_1^\alpha + b z_2^\alpha) = 0$ for any constants $a$ and $b$. In particular, therefore, $\Re(z_1^\alpha + z_2^\alpha) = \Re(e^{\alpha x} + (\ln x)^\alpha) = 0$. We may exploit this fact to quickly compute a resolvent over the smaller differential subfield $\mathbb{Q}(x) = \mathbb{Q}<x> \subset \mathbb{F}$. In general, when $\alpha$ is *not* an integer, there does *not* exist an $\alpha$-resolvent *all of whose terms lie in $\tilde{\mathbb{F}}(\alpha)$ for some smaller differential subfield $\tilde{\mathbb{F}}$* of the polynomial $P_1(t)$'s coefficient differential field $\mathbb{F}$. When $\alpha$ *is* an integer, there often *does* exist an $\alpha$-resolvent over a smaller subfield $\tilde{\mathbb{F}}$, but whose order increases with the absolute magnitude of $\alpha$. In general, the trade-off of desiring a resolvent all of whose terms lie in $\tilde{\mathbb{F}} \subseteq \mathbb{F}, \tilde{\mathbb{F}} \neq \mathbb{F}$ is a resolvent of higher order.

For example, let $y \equiv e^{\alpha \cdot x} + (\ln x)^\alpha$. For each $m \in \mathbb{Z}^+$,

$$D^m y = \alpha^m e^{\alpha \cdot x} + \sum_{k=1}^{m} B_{m,k} \cdot (\alpha)_k \cdot (\ln x)^{\alpha-k} \text{ where} \qquad (1.1)$$

$B_{m,k} = B_{m,k}(x^{-1}, -x^{-2}, 2x^{-3}, \ldots, (-1)^{j-1}(j-1)! \cdot x^{-j}, \ldots (-1)^{k-1}(k-1)! \cdot x^{-k})$ is the Bell polynomial in its arguments, which are negative powers of $x$. The definition of Bell polynomials is given on page 31 of [4]. The falling factorial or Pochhammer symbol $(\alpha)_k$ is defined as $\prod_{i=0}^{k-1}(\alpha - i)$ for $k \in \mathbb{Z}^+$ and 1 for $k = 0$. One can prove that

$B_{m,k}(x^{-1}, -x^{-2}, 2x^{-3}, \ldots, (-1)^{j-1}(j-1)! \cdot x^{-j}, \ldots (-1)^{k-1}(k-1)! \cdot x^{-k})$
$= x^{-m} \cdot B_{m,k}(1, -1, 2, \ldots, (-1)^{j-1}(j-1)!, \ldots (-1)^{k-1}(k-1)!)$.

Define $b_{m,k} \equiv B_{m,k}(1, -1, 2, \ldots, (-1)^{j-1}(j-1)!, \ldots (-1)^{k-1}(k-1)!)$ for $k, m \in \mathbb{Z}^+$ and $k \in [m]$, $b_{0,0} \equiv 1$ and $b_{m,k} \equiv 0$ for all other pairs of values of $m$ and $k$. Then

$$D^m y = \alpha^m e^{\alpha \cdot x} + x^{-m} \cdot \sum_{k=0}^{m} b_{m,k} \cdot (\alpha)_k \cdot (\ln x)^{\alpha-k} \quad \forall m \in \mathbb{Z}_0^+ \qquad (1.2)$$

When $\alpha \in \mathbb{Z}_0^+$, we may use linear algebra to eliminate $e^{\alpha \cdot x}$, $\ln x$, $(\ln x)^2$, $\cdots$, $(\ln x)^\alpha$ from the system of equations (1.2). We get

$$\det \begin{bmatrix} y & 1 & b_{0,0} = 1 & b_{0,1} = 0 & 0 & \cdots & 0 \\ Dy & \alpha & b_{1,0} \cdot x^{-1} = 0 & b_{1,1} \cdot x^{-1} = \alpha \cdot x^{-1} & 0 & \cdots & 0 \\ D^2 y & \alpha^2 & b_{2,0} \cdot x^{-2} = 0 & b_{2,1} \cdot x^{-2} = -\alpha \cdot x^{-2} & b_{2,2} \cdot x^{-2} = \alpha(\alpha-1) \cdot x^{-2} & \cdots & 0 \\ \vdots & \vdots & \vdots & \vdots & \vdots & \ddots & \vdots \\ D^{\alpha+1} y & \alpha^{\alpha+1} & b_{\alpha+1,0} \cdot x^{-\alpha-1} = 0 & b_{\alpha+1,1} \cdot x^{-\alpha-1} & b_{\alpha+1,2} \cdot x^{-\alpha-1} & \cdots & b_{\alpha+1,\alpha+1} \cdot x^{-\alpha-1} \end{bmatrix} = 0$$

(1.3)

Expand the determinant in (1.3), clear denominators, and express the differential equation (1.3) as a linear differential operator $\Re$ acting on $y$ such that $\Re y = 0$. Then $\Re$ is a differential $\alpha$-resolvent of $P_1(t)$ over the differential ring $\mathbb{Z}\{x,\alpha\} = \mathbb{Z}[x]$ (since $\alpha \in \mathbb{Z}_0^+$).
The $j$-th term of $\mathfrak{R}$, by which we mean the coefficient of $D^j y$ in $\mathfrak{R}$, is the $(-1)^{j-1} \cdot (1, j)$-cofactor of the matrix in (1.3).

Observe that the order of $\mathfrak{R}$ depends upon $\alpha$ in (1.3). The order equals $\alpha + 1$.

**Example 1.2.** When $\alpha \notin \mathbb{Z}_0^+$, it is almost certain that there exists no differential $\alpha$-resolvent of $P_1(t)$ whose terms all lie in $\tilde{\mathbb{F}}(\alpha)$ for some proper differential subfield of $\mathbb{F}$. However, the author knows of no proof of this. However, there does exist a second-order differential $\alpha$-resolvent of $P_1(t)$, since $P_1(t)$ is quadratic, which equals

$$\mathfrak{J} \equiv F_{0,2} \cdot D^2 + (F_{0,1} + F_{1,1} \cdot \alpha) \cdot D + (F_{1,0} \cdot \alpha + F_{2,0} \cdot \alpha^2) \quad \mathfrak{J}(e^{\alpha \cdot x}) = 0, \mathfrak{J}((\ln x)^\alpha) = 0 \quad (1.4)$$

One possible formula for a *coefficient-function* -- the individual coefficient of $\alpha^i \cdot D^m$ in $\mathfrak{J}$ -- is the author's powersum formula. This formula is given by

$$F_{0,2} = \det \begin{bmatrix} Dp_1 & 1 \cdot Dp_1 & 1 \cdot p_1 & 1^2 \cdot p_1 \\ Dp_2 & 2 \cdot Dp_2 & 2 \cdot p_2 & 2^2 \cdot p_2 \\ Dp_3 & 3 \cdot Dp_3 & 3 \cdot p_3 & 3^2 \cdot p_3 \\ Dp_4 & 4 \cdot Dp_4 & 4 \cdot p_4 & 4^2 \cdot p_4 \end{bmatrix}, \quad F_{0,1} = \det \begin{bmatrix} D^2 p_1 & 1 \cdot Dp_1 & 1 \cdot p_1 & 1^2 \cdot p_1 \\ D^2 p_2 & 2 \cdot Dp_2 & 2 \cdot p_2 & 2^2 \cdot p_2 \\ D^2 p_3 & 3 \cdot Dp_3 & 3 \cdot p_3 & 3^2 \cdot p_3 \\ D^2 p_4 & 4 \cdot Dp_4 & 4 \cdot p_4 & 4^2 \cdot p_4 \end{bmatrix},$$

and similar expressions for $F_{1,1}$, $F_{1,0}$, and $F_{2,0}$, where $p_q \equiv e^{q \cdot x} + (\ln x)^q$ is the $q$-th powersum of the roots, $e^x$ and $\ln x$, of $P_1(t)$.

**Example 1.3.** Examples 1.1 and 1.2 give us resolvents which annihilate linear combinations over constants of $e^{\alpha \cdot x}$ and $(\ln x)^\alpha$. But, suppose one wants to annihilate arbitrary powers of linear combinations of $e^x$ and $\ln x$ over integers. Let $a$ and $b$ be integers. What does a resolvent for $(a \cdot e^x + b \cdot \ln x)^\alpha$ look like? How would it relate to the resolvent $\mathfrak{R}$ in (1.3) or $\mathfrak{J}$ in (1.4)?

More generally, let $P_1(t) \equiv \prod_{i=1}^{m}(t - z_i) = \sum_{j=0}^{m} \hat{e}_{m-j} \cdot (-1)^{m-j} \cdot t^j$ be a polynomial such that the transcendence degree of the ordinary, non-differential field $\mathbb{Q}(z_1, ..., z_m)$ over $\mathbb{Q}$ is $m$. Let $P_2(t) \equiv \prod_{i=m+1}^{m+n}(t - z_i) = \sum_{j=0}^{n} \breve{e}_{n-j} \cdot (-1)^{n-j} \cdot t^j$ be such that each root $z_i$ satisfies $z_i \in \mathbb{Q}(z_1, ..., z_m), \forall i \in [m+n]$. Let $P_3(t) \equiv \prod_{i=1}^{m+n}(t - z_i) = \sum_{j=0}^{m+n} e_{m+n-j} \cdot (-1)^{m+n-j} \cdot t^j$. The splitting field of $P_3(t)$, $\mathbb{Q}(z_1, ..., z_{m+n}) = \mathbb{Q}(z_1, ..., z_m)$, has infinitely many subfields of finite index. Four of these subfields are $\mathbb{Q}(e_1, ..., e_{m+n})$, $\mathbb{Q}(\hat{e}_1, ..., \hat{e}_m)$, $\mathbb{Q}(\breve{e}_1, ..., \breve{e}_n)$, and $\mathbb{Q}(\hat{e}_1, ..., \hat{e}_m, \breve{e}_1, ..., \breve{e}_n)$ which bear the inclusions
$\mathbb{Q}(\hat{e}_1, ..., \hat{e}_m) \subset \mathbb{Q}(\hat{e}_1, ..., \hat{e}_m, \breve{e}_1, ..., \breve{e}_n) \subset \mathbb{Q}(z_1, ..., z_m)$ and

$\mathbb{Q}(\breve{e}_1,...,\breve{e}_n) \subset \mathbb{Q}(\hat{e}_1,...,\hat{e}_m,\breve{e}_1,...,\breve{e}_n) \subset \mathbb{Q}(z_1,...,z_m)$ and $\mathbb{Q}(e_1,...,e_{m+n}) \subset \mathbb{Q}(\hat{e}_1,...,\hat{e}_m,\breve{e}_1,...,\breve{e}_n)$
Hence, the smallest differential subfields generated by these subfields bear the same inclusions

$\mathbb{Q}<\hat{e}_1,...,\hat{e}_m> \subset \mathbb{Q}<\hat{e}_1,...,\hat{e}_m,\breve{e}_1,...,\breve{e}_n> \subset \mathbb{Q}<z_1,...,z_m>$ and
$\mathbb{Q}<\breve{e}_1,...,\breve{e}_n> \subset \mathbb{Q}<\hat{e}_1,...,\hat{e}_m,\breve{e}_1,...,\breve{e}_n> \subset \mathbb{Q}<z_1,...,z_m>$ and
$\mathbb{Q}<e_1,...,e_{m+n}> \subset \mathbb{Q}<\hat{e}_1,...,\hat{e}_m,\breve{e}_1,...,\breve{e}_n>$. A *joint* resolvent of $P_1(t)$ and $P_2(t)$ is a linear operator whose coefficient-functions are allowed to lie in a larger differential subfield $\mathbb{Q}<\hat{e}_1,...,\hat{e}_m,\breve{e}_1,...,\breve{e}_n>$ than a resolvent of the single polynomial $P_3(t)$, whose coefficient-functions must lie in $\mathbb{Q}<e_1,...,e_{m+n}>$ by definition.

The author tried to prove equality $\mathbb{Q}<e_1,...,e_{m+n}> = \mathbb{Q}<\hat{e}_1,...,\hat{e}_m,\breve{e}_1,...,\breve{e}_n>$ in the special case when $z_i$ is a linear combination of $z_1$ through $z_m$ for each $i \in [m+n]$, but apparently this conjecture is false.

So, let $P_2(t) \equiv \prod_{i=3}^{n+2}(t-(a_i \cdot e^x + b_i \cdot \ln x)) = \sum_{j=0}^{n} \breve{e}_{n-j} \cdot (-1)^{n-j} \cdot t^j$ with
$a_i \cdot b_j - a_j \cdot b_i \neq 0, \forall i \neq j$ and $a_i, b_i \in \mathbb{Q}, \forall i \in [n+2], i \geq 3$. The *coefficient-functions* of a *joint* differential $\alpha$-resolvent of $P_1(t)$ and $P_2(t)$ must lie in
$\mathbb{Q}<e^x + \ln x, e^x \cdot \ln x, \breve{e}_1,...,\breve{e}_n>$. The *terms* of such an $\alpha$-resolvent would lie in $\mathbb{Q}<e^x + \ln x, e^x \cdot \ln x, \breve{e}_1,...,\breve{e}_n>(\alpha)$. A joint $\alpha$-resolvent $\Gamma$ of $P_1(t)$ and $P_2(t)$ would satisfy $\Gamma z^\alpha = 0$ for all roots $z$ of $P_1(t)$ and $P_2(t)$.

Observe that the product polynomial $P_3(t) \equiv P_1(t) \cdot P_2(t)$ has the form
$P_3(t) = \prod_{i=1}^{n+2}(t-(a_i \cdot e^x + b_i \cdot \ln x)) = \sum_{j=0}^{n+2} e_{n+2-j} \cdot (-1)^{n+2-j} \cdot t^j$ with
$a_{n+1}=1, b_{n+1}=0, a_{n+2}=0, b_{n+2}=1$. The coefficient-functions of an $\alpha$-resolvent of $P_3(t)$ must lie in $\mathbb{Q}<e_1,...,e_{n+2}>$. The terms of such an $\alpha$-resolvent would lie in $\mathbb{Q}<e_1,...,e_{n+2}>(\alpha)$.

**Generalization 2 of 3.** Secondly, we will generalize Theorem 8.3 of [6] from resolvents for simple powers of a root $z^\alpha$ to resolvents for general pseudopolynomials involving multiple indeterminate powers $\alpha_{i,j}$. Sketches of this generalization were first presented by the author in a poster at the 14th annual ECCAD conference in [9].

**Example 1.4.** Reconsider Example 1.1. However, this time, we want to find a LODO (Definition 2.5) over $\mathbb{F}(\alpha,\beta) = \mathbb{Q}<e^x + \ln x, e^x \cdot \ln x>(\alpha,\beta)$ which annihilates $e^{\alpha \cdot x} + (\ln x)^\beta$ for any $\alpha$ and $\beta$ where in general $\alpha \neq \beta$, such as when $\alpha$ is an indeterminate variable and $\beta$ is a transcendental number.

There exists a classical way we can compute such a resolvent. Consider the $\alpha$-resolvent $\Im$ of $P_1(t)$ given by (1.4). We need to emphasize that $\Im$ depends upon $\alpha$. So we will write $\Im = \Im_\alpha$. So $\Im_\alpha(e^{\alpha \cdot x}) = 0$. Similarly, $\Im_\beta((\ln x)^\beta) = 0$ if we simply replace $\alpha$ with $\beta$ in (1.4). *We need to use $\Im_\alpha$ and $\Im_\beta$ because their coefficient-functions lie in the required differential field,* $\mathbb{Q} < e^x + \ln x, e^x \cdot \ln x >$. Then the differential consequences, the operators $D\Im_\alpha$, $D^2\Im_\alpha$, $D^3\Im_\alpha$, $D\Im_\beta$, $D^2\Im_\beta$, and $D^3\Im_\beta$, are also differential resolvents of $P_1(t)$. Since $\Im_\alpha$ and $\Im_\beta$ are second-order by (1.4), then $D^3\Im_\alpha$ and $D^3\Im_\beta$ are fifth-order.

When we apply the eight resolvent operators
$\{\Im_\alpha, D\Im_\alpha, D^2\Im_\alpha, D^3\Im_\alpha, \Im_\beta, D\Im_\beta, D^2\Im_\beta, D^3\Im_\beta\}$ to $z_1^\alpha + z_2^\beta = e^{\alpha \cdot x} + (\ln x)^\beta$, we get eight linear combinations of the seven variables $e^{\alpha \cdot x}$, $(\ln x)^\beta$, $(\ln x)^{\beta-1}$, $(\ln x)^{\beta-2}$, $(\ln x)^{\beta-3}$, $(\ln x)^{\beta-4}$, $(\ln x)^{\beta-5}$ over the differential field $\mathbb{Q} < e^x + \ln x, e^x \cdot \ln x > (\alpha, \beta)$. These eight linear combinations equal 0. Hence, there exists some nontrivial linear combination $\Lambda = \Lambda_{\alpha, \beta}$ over $\mathbb{Q} < e^x + \ln x, e^x \cdot \ln x > (\alpha, \beta)$ of the eight resolvents $\{\Im_\alpha, D\Im_\alpha, D^2\Im_\alpha, D^3\Im_\alpha, \Im_\beta, D\Im_\beta, D^2\Im_\beta, D^3\Im_\beta\}$, which annihilates $z_1^\alpha + z_2^\beta = e^{\alpha \cdot x} + (\ln x)^\beta$.

Furthermore, let $\sigma$ be the transposition which switches the roots $z_1 = e^x$ and $z_2 = \ln x$ of $P_1(t)$. By definition of $\Im_\alpha$ as a resolvent of $P_1(t)$, $\Im_\alpha \sigma$ is also a resolvent of $P_1(t)$, that is, $(\Im_\alpha \sigma)e^{\alpha \cdot x} = \Im_\alpha(\sigma e^{\alpha \cdot x}) = \Im_\alpha((\ln x)^\alpha) = 0$. Similarly, $\{\Im_\alpha \sigma, D\Im_\alpha \sigma, D^2\Im_\alpha \sigma, D^3\Im_\alpha \sigma, \Im_\beta \sigma, D\Im_\beta \sigma, D^2\Im_\beta \sigma, D^3\Im_\beta \sigma\}$ are all resolvents of $P_1(t)$. Hence $\Lambda \sigma = \Lambda_{\alpha, \beta} \sigma$ is also a resolvent of $P_1(t)$. Therefore $\Lambda_{\alpha, \beta}(z_1^\beta + z_2^\alpha) = \Lambda_{\alpha, \beta}(e^{\beta \cdot x} + (\ln x)^\alpha) = 0$. Hence, $\Lambda$ is truly a resolvent of $P_1(t)$ for the pseudopolynomial $y = z_1^\alpha + z_2^\beta$.

The method laid out in Example 1.4 is not efficient for computing the terms in a differential resolvent of $P_1(t)$ for the pseudopolynomial $y = z_1^\alpha + z_2^\beta$. However, this classical method *is* useful (in fact, probably necessary) to determine the *form* of such a resolvent, in other words, to determine the number of and powers of $\alpha$ and $\beta$ which would appear in such a resolvent.

To get an idea of the complexity of $\Lambda$ in Example 1.4, start with $\Im$ given by (1.4) in Example 1.2, which has 5 coefficient-functions. Then, one can determine that $D\Im$, $D^2\Im$, and $D^3\Im$ have 9, 12, and 15 coefficient-functions, respectively. However, $\deg_\alpha D^m\Im_\alpha = \deg_\beta D^m\Im_\beta = 2, \forall m \in \mathbb{Z}_0^m$. The order of $D^m\Im_\alpha$ and $D^m\Im_\beta$ is $m + 2$. So $D^m\Im_\alpha$ applied to $e^{\alpha \cdot x} + (\ln x)^\beta$ produces $\beta^i, \forall i \in [m]_0$, and no additional powers of $\alpha$, since $D^m\Im_\alpha$ annihilates the $e^{\alpha \cdot x}$ part. Similarly, $D^m\Im_\beta$ applied to $e^{\alpha \cdot x} + (\ln x)^\beta$ produces $\alpha^i, \forall i \in [m]_0$ and no additional powers of $\beta$, since it annihilates the $(\ln x)^\beta$ part.

Therefore, the eight resolvent operators
$\{\Im_\alpha, D\Im_\alpha, D^2\Im_\alpha, D^3\Im_\alpha, \Im_\beta, D\Im_\beta, D^2\Im_\beta, D^3\Im_\beta\}$, when applied to $e^{\alpha \cdot x} + (\ln x)^\beta$, produce linear combinations of $e^{\alpha \cdot x}$, $(\ln x)^\beta$, $(\ln x)^{\beta-1}$, $(\ln x)^{\beta-2}$, $(\ln x)^{\beta-3}$, $(\ln x)^{\beta-4}$, $(\ln x)^{\beta-5}$
over $\mathbb{Q} < e^x + \ln x, e^x \cdot \ln x > [\alpha, \beta]$. The monomials $\alpha^i \cdot \beta^j$ will range over

$i \in [2]_0$ and $j \in [2]_0$ for $\Im_\alpha(e^{\alpha \cdot x} + (\ln x)^\beta)$

$i \in [2]_0$ and $j \in [3]_0$ for $D\Im_\alpha(e^{\alpha \cdot x} + (\ln x)^\beta)$

$i \in [2]_0$ and $j \in [4]_0$ for $D^2\Im_\alpha(e^{\alpha \cdot x} + (\ln x)^\beta)$

$i \in [2]_0$ and $j \in [5]_0$ for $D^3\Im_\alpha(e^{\alpha \cdot x} + (\ln x)^\beta)$

$i \in [2]_0$ and $j \in [2]_0$ for $\Im_\beta(e^{\alpha \cdot x} + (\ln x)^\beta)$

$i \in [3]_0$ and $j \in [2]_0$ for $D\Im_\beta(e^{\alpha \cdot x} + (\ln x)^\beta)$

$i \in [4]_0$ and $j \in [2]_0$ for $D\Im_\beta(e^{\alpha \cdot x} + (\ln x)^\beta)$

$i \in [5]_0$ and $j \in [2]_0$ for $D^2\Im_\beta(e^{\alpha \cdot x} + (\ln x)^\beta)$

Eliminating the $e^{\alpha \cdot x}$ and $(\ln x)^{\beta-m}$ terms is the equivalent of computing the determinant of an $8 \times 8$ matrix with entries in $\mathbb{Q} < e^x + \ln x, e^x \cdot \ln x > [\alpha, \beta]$. The final result is that $\Lambda$ is a fifth-order resolvent whose terms lie in $\mathbb{Q} < e^x + \ln x, e^x \cdot \ln x > [\alpha, \beta]$ with all monomials $\alpha^i \cdot \beta^j$ appearing for each $i, j \in [2+2+2+2+3+4+5]_0 = [20]_0$. Thus, each term of $\Lambda$ will have $21 \cdot 21 = 441$ coefficient functions. In actuality, the fifth-order and zeroeth order terms will most likely have fewer. Only a direct computation will prove that. So, we estimate that the six-term resolvent $\Lambda$ will have about $6 \cdot 441 = 2646$ coefficient-functions.

**Generalization 3 of 3.** In [8] the author generalized the definition of a resolvent of a single polynomial to include nonzero differential equations with *nonconsecutive* derivatives. We will include this generalization in this article.

Then we will demonstrate in Theorem 4.2 that the author's powersum formula [5] can be used to compute joint resolvents when they exist. Finally, in Section 5, we will compute an elementary joint resolvent using Theorem 4.2.

We will *not* attempt to prove rigorously that the powersum formula can always be used to produce a nonzero differential equation. At present, this is a conjecture. (The definition of a differential resolvent contains the condition that the differential equation be *nonzero*.) No attempt will be made to determine the specific *form* (the specific sequence of powers of $\alpha_{i,j}$ in a resolvent). Example 1.5 will demonstrate that the powersum formula can yield an identically zero LODO (Definition 2.5).

**Example 1.5. Differential field with nonzero characteristic**

Let $\mathbb{F}_3$ be a differential field of characteristic 3. Let $x \in \mathbb{F}_3$. Let $P(t) \equiv t^3 + x \cdot t - 1$. Let $P(z) = 0$. So $z^3 + x \cdot z - 1 = 0$. So $(3z^2 + x) \cdot Dz + z = 0$. Since $3z^2 \equiv 0 \bmod 3$, then $x \cdot Dz + z = 0$. Let $y \equiv z^\alpha$. So $x \cdot Dy + \alpha \cdot y = 0$. So $x \cdot D + \alpha$ is a first-order linear differential $\alpha$-resolvent of $P(t)$.

In order to recover this same resolvent using the powersum formula, set up undetermined coefficient-functions $r_1$ and $r_0$ with $r_1, r_0 \in \mathbb{F}_3$ such that

$r_1 \cdot Dy + \alpha \cdot r_0 \cdot y = 0$. The first three elementary symmetric functions of the roots of $P(t)$ are $e_1 = 0$, $e_2 = x$, $e_3 = 1$. Using the formula given on page 28 of Macdonald [4], one can determine that the first three powersums of the roots are $p_1 = 0$, $p_2 = -2x$, and $p_3 = 3 \equiv 0 \bmod 3$.

The powersum formula for the resolvent in Example 1.5 is $r_1 = \alpha \cdot p_\alpha$ and $r_0 = -Dp_\alpha$ where we specialize $\alpha$ to an integer. Observe that since the resolvent has only 2 coefficient-functions, the powersum formula is expressed as the determinant of a 1 x 1 matrix – a single entry. If we specialize $\alpha \to 1$, then $r_1 = 1 \cdot p_1 = 0$ and $r_0 = -Dp_1 = 0$. If we specialize $\alpha \to 3 \equiv 0 \bmod 3$, then $r_1 = 3 \cdot p_3 \equiv 0$ and $r_0 = -Dp_3 = 0$. So, the powersum formula yields an identically zero differential operator for these two choices of $\alpha$. However, if we specialize $\alpha \to 2$, then $r_1 = 2 \cdot p_2 = -4x \equiv 2x \bmod 3$ and $r_0 = -Dp_2 = 2$. Now divide $r_1$ and $r_0$ by the common factor of 2 to recover the resolvent $x \cdot D + \alpha$.

Example 1.4 has demonstrated that determining the specific form of resolvents of given polynomials is an enormous task. It has been done only in the simplest of cases, such as in Theorem 8.3 of [6]. Theorems 1, 2 and 3 of [1] are collectively a tremendous attack on this problem for resolvents of a single polynomial $P(X,Y) = 0$, polynomial in the independent variable $X$ as well as $Y$.

**Example 1.6.** Let $P(t) \in \mathbb{Z}[t, x]$ have non-constant roots. There can exist no $\alpha$-resolvent of $P(t)$ of the form $f_1 \cdot Dz^\alpha + (f_2 + f_3 \cdot \alpha^m) \cdot z^\alpha = 0$ for any $m \neq 1$. The proof is easy. Simply integrate $\dfrac{Dz^\alpha}{z^\alpha} = \dfrac{f_2 + f_3 \cdot \alpha^m}{f_1}$. We would get $z^\alpha = c^\alpha \cdot h(x) \cdot (r(x))^{\alpha^m}$ for some functions $h(x)$ and $r(x)$ of $x$ and constant $c \equiv z|_{x=0}$. Since $\alpha$ is independent of $x$, such a functional relationship is impossible. Otherwise, we would need $(r(x))^{\alpha^m} = (r(x))^\alpha, \forall x$. This is possible only if $r(x) \equiv 0, \forall x$ or $r(x) \equiv 1, \forall x$. Similarly, $h(x) \equiv 1, \forall x$ would be necessary. Thus, $z = z(x)$ would have to be constant with respect to $x$, contradicting our hypothesis that $P(t)$ has non-constant roots.

No attempt will be made in this paper to state a general theorem to help us factor each of the terms (more specifically, the coefficient-functions) of the joint resolvent. In [8] the author demonstrated how the main factorization theorem, Theorem 8.2 page 3089, in [7], can be generalized to resolvents with nonconsecutive derivatives. Work is being done to generalize this factorization theorem to include joint resolvents for pseudopolynomials involving multiple indeterminate powers $\alpha_{i,j}$. Such a formal matrix factorization will be crucial in order to make the computation of joint resolvents possible. Currently, unfortunately, the computation of joint resolvents for pseudopolynomials is not computable on personal computers except for only the simplest of polynomials, such as the linear polynomials in Example 6.1.

Section 5 will discuss other lines of research based on this work, such as determining the transcendentality of certain numbers expressed as solutions of a differential resolvent evaluated at certains points.

**2. Basic definitions concerning resolvents**

Let $\mathbb{F}$ be an ordinary differential field of characteristic 0 with derivation $D$. Let $\{\alpha_i\}_{i=1}^{i=M}$ be a collection of $M$ distinct transcendental constants (with respect to $D$) over $\mathbb{F}$ such that $tr\deg_{\mathbb{F}} \mathbb{F}(\{\alpha_i\}_{i=1}^{i=M}) = M$. So $D\alpha_i = 0, \forall i \in [M]$.

**Definition 2.2.** Define $\sum_{i=1}^{i=M} a_i \cdot t_i^{\alpha_i}$ with $a_i \in \mathbb{F}$ to be a *multivariate pseudo-polynomial in the variables $t_i$ over $\mathbb{F}$*.

Let $\{e_{k,j}\}_{k \in [n_j]_0, j \in [L]} \subset \mathbb{F}$ be a collection of $\sum_{j=1}^{L} n_j$ not necessarily distinct elements of $\mathbb{F}$ with $e_{0,j} \equiv 1, \forall j \in [L]$. These elements will be the coefficients of monic univariate polynomials $P_j$. There is no critical reason now to take the polynomials $P_j$ to be monic. That need arises when we examine deeper lines of research when we use the powersum formula Theorem 4.2 to prove the smallest *rings* in which differential resolvents can lie. Hence, it is safest now to remain consistent with standard notation and practice.

**Definition 2.3.** Define $P_j(t) \equiv \sum_{k=0}^{n_j} (-1)^{n_j - k} \cdot e_{n_j - k, j} \cdot t^k$.

Let $u_j$ represent any one of the roots of $P_j$. So $P_j(u_j) = 0$. So $e_{k,j}$ equals the $k$-th elementary symmetric function of all $n_j$ roots of $P_j$.

Define $\Lambda \equiv \{e_{k,j}\}_{k \in [n_j], j \in [L]}$ to be the set of all $\sum_{j=1}^{L}(n_j + 1)$ coefficients of all the $P_j$.

Define $\mathbb{Q}\{\Lambda\}$ to be the differential ring generated by the rational numbers and $\Lambda$.

**Definition 2.4.** For the remainder of this paper we will assume $\mathbb{F} = \mathbb{Q} < \Lambda >$. In other words, $\mathbb{F}$ will be the finitely differentially generated field of characteristic 0 generated by the coefficients of all the polynomials $\{P_j\}_{j \in [L]}$. This means that $\mathbb{F}$ will be the *smallest differential field* containing the coefficients of the $\{P_j\}_{j \in [L]}$.

**Definition 2.5.** We will use the acronym *LODO* for linear ordinary differential operator of finite order.

**Definition 2.6.** We say a LODO, $\sum_{m \geq 0, m \in S \subset \mathbb{Z}_0^+} r_m \cdot D^m$, is *zero* or *identically zero,* if $r_m = 0, \forall m \in S \subset \mathbb{Z}_0^+$, for whatever indexing subset $S$ of the nonnegative integers is chosen.

Definition 2.6 for a zero LODO is very different from saying $\sum_{m \in S \subset \mathbb{Z}_0^+} r_m \cdot D^m y = 0$ for some solution $y$ of the LODO $\sum_{m \in S \subset \mathbb{Z}_0^+} r_m \cdot D^m$.

Theorem 8.3 page 2876 of [6] and Theorem 37 page 67 of [10] assert that there exists a LODO $\sum_{m=0}^{n} \sum_{j=0}^{\frac{n(n-1)}{2}+1-m} r_{j,m} \cdot \alpha^j \cdot D^m$, with all $r_{j,m} \in \mathbb{F}$, not all 0, satisfied by $z^\alpha$ for all $n$ roots $z$ of a *single polynomial* $P$ for a *single* indeterminate constant $\alpha$.

**Definition 2.7.** Let $\alpha$ be a constant lying in a differential field extension of $\mathbb{F}$. A LODO $\sum_{m \in S \subset \mathbb{Z}_0^+} r_m \cdot D^m$ associated with a polynomial $P$ is called an $\alpha$-*resolvent*, or a *differential resolvent*, or a *linear differential resolvent* of $P$ if all $r_m \in \mathbb{F}$, the $r_m$ are not all 0, and the LODO annihilates all $n$ roots $z$ of $P$ when raised to the indeterminate constant power $\alpha$. In other words, $\sum_{m \in S \subset \mathbb{Z}_0^+} r_m \cdot D^m z^\alpha = 0$ for all $n$ roots $z$ of $P$.

The definition of a differential resolvent given by Definition 2.7 is unrelated to the notion of a "resolvent of a differential prime ideal", supposedly also called a "differential resolvent", by Joseph Fels Ritt in his book *Differential Algebra*.

We will typically denote differential resolvents by capital letters $R$ or $\mathfrak{R}$ and affix various subscripts and superscripts to them as necessary. We will typically denote the coefficient of $D^m$ by small letter $r_m$.

**Definition 2.8.** We will refer to $r_m$ as a *term* of the resolvent.

If we expand the resolvent as a polynomial in $\alpha$, then we will typically denote the coefficient of $\alpha^j D^m$ by $r_{j,m}$.

**Definition 2.9.** We will refer to $r_{j,m}$ as a *coefficient-function* of the resolvent.

When we write a resolvent, or any LODO, in the form $\sum_{m \in S \subset \mathbb{Z}_0^+} r_m \cdot D^m$, it is possible that some of the $r_m$ are 0. If we wish to notate the LODO showing only the nonzero terms, we will write the LODO in the form $\sum_{k \in [n]} r_{m_k} \cdot D^{m_k} = 0$ for some $n \in \mathbb{Z}^+$ where $n = |S|$. This way, we can take the subindex $k$ on the derivatives $m_k$ to run over consecutive values $[n]$ even if the $m_k$ are not consecutive.

The first purpose of this paper is to extend the definition of a resolvent from Definition 2.7. More generally, let $\sum_{j=1}^{j=M} a_j \cdot \prod_{i=1}^{L} t_i^{\alpha_{i,j}}$ be a *multivariate pseudopolynomial* over $\mathbb{F}$, that is, $a_j \in \mathbb{F}$. As will be explained in the proof of Theorem 3.1, it is sufficient to consider pseudopolynomials of the form $\sum_{j=1}^{j=M}\prod_{i=1}^{L} t_i^{\alpha_{i,j}}$ without the $a_j$. However, leaving the $a_j$ in sometimes helps remind the reader that this pseudopolynomial is considered to be a polynomial over the differential field $\mathbb{F}$.

Let $\{u_i\}_{i=1}^{i=L}$ be a set of $L$ algebraic elements over $\mathbb{F}$, generally with relations among them. Let $P_i(t)$ be the minimal polynomial over $\mathbb{F}$ of $u_i$. Let $\sigma_i$ represent any element of the Galois group of $P_i(t)$. So $\sigma_i$ sends $u_i$ to some other root $u_i{}'$ of $P_i(t)$.

**Definition 2.10.** Let $\sigma$ denote a *combined permutation* of all the roots of all the polynomials $P_i(t)$. Then $\sigma$ is defined by $\sigma\left(\{u_i\}_{i=1}^{i=L}\right) = \{u_i{}'\}_{i=1}^{i=L}$. In other words, $\sigma$ maps the roots of $P_i(t)$ to the roots of $P_i(t)$, not to $P_j(t)$ if $j \neq i$.

**3. Existence of joint resolvents of multiple polynomials for a pseudopolynomial**

In applications of Theorem 3.1, one is more likely to start with a collection of univariate polynomials $\{P_i(t)\}_{i=1}^{L}$ rather than with a common differential coefficient field $\mathbb{F}$. Thereofre, we will implicitly assume that, if a collection of univariate polynomials $\{P_i(t)\}_{i=1}^{L}$ has a collection of individual coefficient fields $\{\mathbb{F}_i(t)\}_{i=1}^{L}$, then all $\{\mathbb{F}_i(t)\}_{i=1}^{L}$ have the same characteristic.

**Theorem 3.1.** Let $\{P_i(t)\}_{i=1}^{L}$ be a set of $L$ univariate polynomials whose coefficients lie in a differential field $\mathbb{F}$ with $\deg_t P_i(t) = d_i$. Let $m_1 < m_2 < \cdots < m_N$ be any strictly increasing sequence of $N$ nonnegative integers where $N \geq \prod_{i=1}^{L} d_i$. Take $\mathbb{F}$ to be the smallest differential field containing (generated by) all the coefficients of the $\{P_i(t)\}_{i=1}^{L}$ and $\mathbb{Q}$. Then there exists a nonzero resolvent $\sum_{k=1}^{N}\sum_{s=1}^{h_k} r_{k,s} \cdot \prod_{i=1}^{L}\prod_{j=1}^{M} \alpha_{i,j}^{l_{i,j,k,s}} \cdot D^{m_k}$ with coefficient-functions $r_{k,s} \in \mathbb{F}^{\#}$ such that $\sum_{k=1}^{N}\sum_{s=1}^{h_k} r_{k,s} \cdot \prod_{i=1}^{L}\prod_{j=1}^{M} \alpha_{i,j}^{l_{i,j,k,s}} \cdot D^{m_k}(\sigma y) = 0$ holds *for* the multipseudopolynomial $y \equiv \sum_{j=1}^{j=M} a_j \cdot \prod_{i=1}^{L} u_i^{\alpha_{i,j}}$ for every combined permutation $\sigma$ (Definition 2.10) of the roots *of* $\{P_i\}_{i=1}^{L}$. In the multipseudopolynomial $y$, we use the notation such that $l_{i,j,k,s} \in \mathbb{Z}_0^{+}, \forall i \in [L], \forall j \in [M], \forall s \in [h_k]$ and $\sum_{i=1}^{L}\sum_{j=1}^{M} l_{i,j,k,s} = \rho$ and $N - 1 \leq \rho \leq \sum_{k'=1}^{N} m_{k'} - m_k$.

**Organization of the proof of Theorem 3.1.** We will differentiate the multipseudopolynomial $y \equiv \sum_{j=1}^{j=M} a_j \cdot \prod_{i=1}^{L} u_i^{\alpha_{i,j}}$ a sufficient number of times, namely, at least $M \cdot \prod_{i=1}^{L} d_i$ times. We will reduce all integer powers of $u_i$ and all derivatives of $u_i$ modulo $P_i(t)$. Hence, we will express all derivatives of $y$ as linear combinations of the monomials $\{\prod_{i=1}^{L} u_i^{\alpha_{i,j}}\}_{j=1}^{M}$ over the joint differential coefficient field $\mathbb{F}$ of the polynomials $\{P_i(t)\}_{i=1}^{L}$ adjoined with all the indeterminates $\{\alpha_{i,j}\}$. Elementary vector space theory will guarantee a nontrivial relation among these derivatives of $y$ over $\mathbb{F}[\alpha]$. This relation – more specifically, the linear differential operator which when applied to $y$ yields zero – will be a differential resolvent of the $\{P_i(t)\}_{i=1}^{L}$ for $y$.

**Proof.** First observe that if $a_j \in \mathbb{F}$ in $y = \sum_{j=1}^{j=M} a_j \cdot \prod_{i=1}^{L} u_i^{\alpha_{i,j}}$, then we may add to our collection of polynomials $\{P_i\}_{i=1}^{L}$ the polynomial $\tilde{P}(t) \equiv a_j - t$ with root $\tilde{u} = a_j$ for each $j \in [M]$. So, it is sufficient to prove Theorem 3.1 for the multipseudopolynomial $y = \sum_{j=1}^{j=M} \prod_{i=1}^{L} u_i^{\alpha_{i,j}}$. We cannot directly take the logarithmic derivative of $y$ in order to prove there exists a resolvent for it. We would immediately run into the problem, for instance, when $y = z^\alpha + z^\beta$, that $\frac{Dy}{y} = \frac{\alpha z^{\alpha-1} + \beta z^{\beta-1}}{z^\alpha + z^\beta} \notin \mathbb{F}[\alpha, \beta, z]$. In spite of this difficulty, there are any number of ways to proceed. We will choose a method which is the most "symmetrical" in the derivatives of $y$.

So, with no loss of generality, consider $y \equiv \sum_{j=1}^{j=M} \prod_{i=1}^{L} u_i^{\alpha_{i,j}}$. (3.1)

Define $v_j \equiv \prod_{i=1}^{L} u_i^{\alpha_{i,j}}$. (3.2)

So $y = \sum_{j=1}^{M} v_j$ and $\frac{Dv_j}{v_j} = \sum_{i=1}^{L} \alpha_{i,j} \cdot \frac{Du_i}{u_i} = \sum_{i=1}^{L} \alpha_{i,j} \cdot \tilde{P}_i(u_i)$ for some polynomial $\tilde{P}_i(t) \in \mathbb{F}[t]$ whose coefficients lie in the subfield of $\mathbb{F}$ differentially generated by the coefficients of $P_i(t)$. We use the $\sim$ symbol over $P_i(t)$ to distinguish the polynomials $\tilde{P}_i(t)$ from $P_i(t)$. So $Dy = \sum_{j=1}^{j=M} Dv_j = \sum_{j=1}^{j=M} v_j \cdot \frac{Dv_j}{v_j} = \sum_{j=1}^{j=M} v_j \cdot \sum_{i=1}^{L} \alpha_{i,j} \cdot \tilde{P}_i(u_i)$. So $\sigma y = \sum_{j=1}^{j=M} \sigma v_j$ and hence $D^m(\sigma y) = \sum_{j=1}^{j=M} D^m(\sigma v_j)$ for every $m \in \mathbb{Z}$ with $m > 0$.

Observe that $Dy$ is first order in the constants $\alpha_{i,j}$. We make the induction hypothesis that the form of the $m$-th derivative of $y$ is

$$D^m y = \sum_{j=1}^{j=M} v_j \cdot \left\{ \begin{array}{l} \sum_{i_1=1}^{L} \alpha_{i_1,j} \cdot \tilde{P}_{1,i_1,m}(u_{i_1}) + \sum_{1 \leq i_1 < i_2 \leq L} \alpha_{i_1,j} \cdot \alpha_{i_2,j} \cdot \tilde{P}_{2,i_1,i_2,m}(u_{i_1}, u_{i_2}) + \cdots \\ + \sum_{1 \leq i_1 < \cdots < i_m \leq L} \left( \prod_{s=1}^{m} \alpha_{i_s,j} \right) \cdot \tilde{P}_{m,i_1,i_2,\cdots,i_m,m}(u_{i_1}, u_{i_2}, \cdots, u_{i_m}) \end{array} \right\}. \quad (3.3)$$

We prove hypothesis (3.3) true by differentiating it

$$D^{m+1} y = \sum_{j=1}^{j=M} Dv_j \cdot \left\{ \begin{array}{l} \sum_{i_1=1}^{L} \alpha_{i_1,j} \cdot \tilde{P}_{1,i_1,m}(u_{i_1}) + \sum_{1 \leq i_1 < i_2 \leq L} \alpha_{i_1,j} \cdot \alpha_{i_2,j} \cdot \tilde{P}_{2,i_1,i_2,m}(u_{i_1}, u_{i_2}) + \cdots \\ + \sum_{1 \leq i_1 < \cdots < i_m \leq L} \left( \prod_{s=1}^{m} \alpha_{i_s,j} \right) \cdot \tilde{P}_{m,i_1,i_2,\cdots,i_m,m}(u_{i_1}, u_{i_2}, \cdots, u_{i_m}) \end{array} \right\}$$

$$+ \sum_{j=1}^{j=M} v_j \cdot \left\{ \begin{array}{l} \sum_{i_1=1}^{L} \alpha_{i_1,j} \cdot D\tilde{P}_{1,i_1,m}(u_{i_1}) + \sum_{1 \leq i_1 < i_2 \leq L} \alpha_{i_1,j} \cdot \alpha_{i_2,j} \cdot D\tilde{P}_{2,i_1,i_2,m}(u_{i_1}, u_{i_2}) + \cdots \\ + \sum_{1 \leq i_1 < \cdots < i_m \leq L} \left( \prod_{s=1}^{m} \alpha_{i_s,j} \right) \cdot D\tilde{P}_{m,i_1,i_2,\cdots,i_m,m}(u_{i_1}, u_{i_2}, \cdots, u_{i_m}) \end{array} \right\} \quad (3.4)$$

where we now replace $Dv_j$ on the inside of the first summation of the first line of (3.3) with $\sum_{i=1}^{L} \alpha_{i,j} \cdot \tilde{P}_i(u_i)$ to get (3.5)

$$D^{m+1} y = \sum_{j=1}^{j=M} v_j \cdot \sum_{i=1}^{L} \alpha_{i,j} \cdot \tilde{P}_i(u_i) \cdot \left\{ \begin{array}{l} \sum_{i_1=1}^{L} \alpha_{i_1,j} \cdot \tilde{P}_{1,i_1,m}(u_{i_1}) + \sum_{1 \leq i_1 < i_2 \leq L} \alpha_{i_1,j} \cdot \alpha_{i_2,j} \cdot \tilde{P}_{2,i_1,i_2,m}(u_{i_1}, u_{i_2}) + \cdots \\ + \sum_{1 \leq i_1 < \cdots < i_m \leq L} \left( \prod_{s=1}^{m} \alpha_{i_s,j} \right) \cdot \tilde{P}_{m,i_1,i_2,\cdots,i_m,m}(u_{i_1}, u_{i_2}, \cdots, u_{i_m}) \end{array} \right\}$$

$$+ \sum_{j=1}^{j=M} v_j \cdot \left\{ \begin{array}{l} \sum_{i_1=1}^{L} \alpha_{i_1,j} \cdot D\tilde{P}_{1,i_1,m}(u_{i_1}) + \sum_{1 \leq i_1 < i_2 \leq L} \alpha_{i_1,j} \cdot \alpha_{i_2,j} \cdot D\tilde{P}_{2,i_1,i_2,m}(u_{i_1}, u_{i_2}) + \cdots \\ + \sum_{1 \leq i_1 < \cdots < i_m \leq L} \left( \prod_{s=1}^{m} \alpha_{i_s,j} \right) \cdot D\tilde{P}_{m,i_1,i_2,\cdots,i_m,m}(u_{i_1}, u_{i_2}, \cdots, u_{i_m}) \end{array} \right\} \quad (3.5)$$

When the inner sum $\sum_{i=1}^{L} \alpha_{i,j} \cdot \tilde{P}_i(u_i)$ in (3.5) gets multiplied by the term in brackets in the first line of (3.5), it yields a term of total degree in the $\alpha_{i,j}$ of $m+1$, because the total degree in the $\alpha_{i,j}$ of the bracketed terms is $m$ and the total degree of $\sum_{i=1}^{L} \alpha_{i,j} \cdot \tilde{P}_i(u_i)$ is 1. Hence, the form of $D^{m+1} y$ is

$$D^{m+1}y = \sum_{j=1}^{j=M} v_j \cdot \left\{ \begin{array}{l} \sum_{i_1=1}^{L} \alpha_{i_1,j} \cdot \tilde{P}_{1,i_1,m+1}(u_{i_1}) + \sum_{1 \leq i_1 < i_2 \leq L} \alpha_{i_1,j} \cdot \alpha_{i_2,j} \cdot \tilde{P}_{2,i_1,i_2,m+1}(u_{i_1},u_{i_2}) + \cdots \\ + \sum_{1 \leq i_1 < \cdots < i_{m+1} \leq L} \left( \prod_{s=1}^{m+1} \alpha_{i_s,j} \right) \cdot \tilde{P}_{m+1,i_1,i_2,\cdots,i_{m+1},m+1}(u_{i_1},u_{i_2},\cdots,u_{i_{m+1}}) \end{array} \right\} \quad (3.6)$$

where $\sum_{i_1=1}^{L} \alpha_{i_1,j} \cdot \tilde{P}_{1,i,m+1}(u_i) = \sum_{i_1=1}^{L} \alpha_{i_1,j} \cdot D\tilde{P}_{1,i_1,m}(u_{i_1})$ and (3.7a)

$$\sum_{1 \leq i_1 < i_2 \leq L} \alpha_{i_1,j} \cdot \alpha_{i_2,j} \cdot \tilde{P}_{2,i_1,i_2,m+1}(u_{i_1},u_{i_2})$$
$$= \left( \sum_{i=1}^{L} \alpha_{i,j} \cdot \tilde{P}_i(u_i) \right) \cdot \left( \sum_{i_1=1}^{L} \alpha_{i_1,j} \cdot \tilde{P}_{1,i,m}(u_i) \right) + \sum_{1 \leq i_1 < i_2 \leq L} \alpha_{i_1,j} \cdot \alpha_{i_2,j} \cdot D\tilde{P}_{2,i_1,i_2,m}(u_{i_1},u_{i_2}) \quad (3.7b)$$

and… omitting the set of equations relating the polynomial $\tilde{P}_{k,i_1,\ldots,i_k,m+1}(u_{i_1},\ldots,u_{i_k})$ to the polynomials $\{\tilde{P}_{j,i_1,\ldots,i_j,m}(u_{i_1},\ldots,u_{i_j})\}_{j=1}^{k}$ for $3 \leq k \leq m$ …

$$\sum_{1 \leq i_1 < \cdots < i_{m+1} \leq L} \left( \prod_{s=1}^{m+1} \alpha_{i_s} \right) \cdot \tilde{P}_{m+1,i_1,i_2,\cdots,i_{m+1},m+1}(u_{i_1},u_{i_2},\cdots,u_{i_{m+1}}) =$$
$$\left( \sum_{i=1}^{L} \alpha_{i,j} \cdot \tilde{P}_i(u_i) \right) \cdot \sum_{1 \leq i_1 < \cdots < i_m \leq L} \left( \prod_{s=1}^{m} \alpha_{i_s,j} \right) \cdot \tilde{P}_{m,i_1,i_2,\cdots,i_m,m}(u_{i_1},u_{i_2},\cdots,u_{i_m}) \quad (3.7n)$$
$$+ \sum_{1 \leq i_1 < \cdots < i_m \leq L} \left( \prod_{s=1}^{m} \alpha_{i_s,j} \right) \cdot D\tilde{P}_{m,i_1,i_2,\cdots,i_m,m}(u_{i_1},u_{i_2},\cdots,u_{i_m})$$

Since equation (3.6) is equation (3.3) with the induction index $m$ replaced by $m+1$, the induction hypothesis is true for all positive integers $m$.

Let $\alpha$ stand collectively for all $L \cdot M$ $\alpha_{i,j}$'s in equations (3.7).

Let E denote the differential field $\mathbb{F}(u_1,\ldots,u_L,\alpha)$. We use the notation E to stand for an *extension* of $\mathbb{F}(\alpha)$. Since $[\mathbb{F}(u_1,\ldots,u_L):\mathbb{F}] \leq \prod_{i=1}^{L} d_i$, it follows that

$$[\mathrm{E}:\mathbb{F}(\alpha)] = [\mathbb{F}(u_1,\ldots,u_L,\alpha):\mathbb{F}(\alpha)] \leq [\mathbb{F}(u_1,\ldots,u_L):\mathbb{F}] \leq \prod_{i=1}^{L} d_i.$$

Let $V$ denote the vector space $\mathrm{E} \cdot v_1 + \mathrm{E} \cdot v_2 + \cdots \mathrm{E} \cdot v_M$ generated by the $\{v_j\}_{j=1}^{M}$ over E. Then $\dim_{\mathrm{E}} V \leq M$. So $\dim_{\mathbb{F}(\alpha)} V = \dim_{\mathrm{E}} V \cdot \dim_{\mathbb{F}(\alpha)} \mathrm{E} \leq M \cdot \prod_{i=1}^{L} d_i$.

Since each $D^m y \in V$ by equation (3.3), we are guaranteed by basic linear algebra that there exists a nontrivial relation among any $N$ derivatives $\{D^{m_k} y\}_{k=1}^{N}$ of $y$ if $N \geq M \cdot \prod_{i=1}^{L} d_i$. Expand the right side of equation (3.3) as sums of monomials of the $v_j$ and integer powers of $u_i$. Eliminate all monomials consisting of $v_j$ (the transcendental powers of $u_i$) and integer powers of $u_i$ with the use of a(n extremely large) determinant.

The cofactor in this determinant of $D^{m_k} y$ will have degree in the $\alpha_{i,j}$'s equal to $\left(\sum_{k'=1}^{N} m_{k'}\right) - m_k$. Since every entry in this determinant has no zeroeth degree term in the $\alpha_{i,j}$'s, there will be no $\rho$-th order terms in $\alpha_{i,j}$'s for $\rho < N-1$. This proves that for each $L \cdot M$ tuple $\vec{l}_{k,s} = (l_{1,1,k,s}, \ldots, l_{L,M,k,s})$ of nonnegative integer powers, there exists a differential resolvent such that $\sum_{i=1}^{L} \sum_{j=1}^{M} l_{i,j,k,s} = \rho$ with $N - 1 \leq \rho \leq \sum_{k'=1}^{N} m_{k'} - m_k$.

The resolvent is written in its "nonzero" form as $\sum_{k=1}^{N} \sum_{s=1}^{h_k} r_{k,s} \cdot \prod_{i=1}^{L} \prod_{j=1}^{M} \alpha_{i,j}^{l_{i,j,k,s}} \cdot D^{m_k}$.

That means all coefficient-functions $r_{k,s} \neq 0$ and that there are $\sum_{k=1}^{N} h_k$ such nonzero coefficient-functions. **Q.E.D.**

**Comment on the proof of Theorem 3.1.** There might be relations among the $\alpha_{i,j}$, such as $\alpha_{1,1} = \sqrt{5}$, $\alpha_{1,2} = -\sqrt{5}$, $\alpha_{2,1} = \pi$, $\alpha_{2,2} = -\pi$. But, in the "worst case" there are no relations among the $\alpha_{i,j}$. We say "worst case", because a joint resolvent for a pseudopolynomial in which no relations exist among the $\alpha_{i,j}$ will have much greater order and many more terms than a joint resolvent for a pseudopolynomial in which algebraic relations exist among the $\alpha_{i,j}$. In this worst case, $\mathrm{tr\,deg}_{\mathbb{F}} \mathbb{F}(\alpha) = L \cdot M$ and $\dim_E V = M$. Also, in the "worst case" there are no relations among the polynomials $\{P_i\}_{i=1}^{L}$, so $\dim_{\mathbb{F}(\alpha)} E = \prod_{i=1}^{L} d_i$. So, if both these worst cases hold, then

$$\dim_{\mathbb{F}(\alpha)} V = M \cdot \prod_{i=1}^{L} d_i.$$

Observe that the resolvent in Theorem 3.1 is itself a multivariate polynomial in the indeterminate constant powers $\alpha_{i,j}$.

**Definition 3.2.** The differential resolvent $\sum_{k=1}^{N} \sum_{s=1}^{h_k} r_{k,s} \cdot \prod_{i=1}^{L} \prod_{j=1}^{M} \alpha_{i,j}^{l_{i,j,k,s}} \cdot D^{m_k}$ of Theorem 3.1 is called a *joint differential resolvent of* the polynomials $\{P_i\}_{i=1}^{L}$ *for* the pseudopolynomial

$$y \equiv \sum_{j=1}^{j=M} a_j \cdot \prod_{i=1}^{L} u_i^{\alpha_{i,j}}.$$

In fact, a joint differential resolvent of two or more polynomials may be a resolvent for a pseudo (or regular) polynomial of the roots of one of the polynomials and a resolvent for a *different* pseudo (or regular) polynomial of the roots of the other polynomial. Example 3.3 explains.

**Example 3.3.** The differential ring $\mathbb{Z}\{x\}$ is the same set as the ordinary polynomial ring $\mathbb{Z}[x]$ if $Dx \equiv 1$. Let $P(t) \equiv \sum_{k=0}^{n} a_k \cdot t^k$ be a univariate polynomial over $\mathbb{Z}[x]$. Assume the greatest common divisor (gcd) of the $a_k$ over $\mathbb{Z}[x]$ is 1. If all the coefficients $a_k$ of a polynomial $P(t) = \sum_{k=0}^{n} a_k \cdot t^k$ over $\mathbb{Z}\{x\}$ were constant, in other words, were all integers, then all the roots of $P(t)$ would be constant, hence, independent of $x$. If at least one of the $a_k$ would involve $x$, in other words, if at least one of the $a_k$ would be a polynomial in $x$ of degree >0, then some derivative of that $a_k$ would be $m_k \cdot x$ for some integer $m_k \in \mathbb{Z}^+$. Hence, if all the coefficients $a_k$ of a polynomial $P(t) = \sum_{k=0}^{n} a_k \cdot t^k$ lie in $\mathbb{Z}[x]$ but are not all constant, then they and all their derivatives generate the differential coefficient ring $m \cdot \mathbb{Z}[x] = m \cdot \mathbb{Z}\{x\}$ where $m \equiv \gcd\{m_k, \forall k \in [n]_0\}$. Regardless of the value of $m$, the differential coefficient *field* of $P$ is $\mathbb{Q}<x> = \mathbb{Q}(x)$. Hence, if all the coefficients of another polynomial $P_2(t)$ lie in $\mathbb{Z}[x] \subset \mathbb{Q}(x)$, then each of the coefficients of $P_2(t)$ is generated by the coefficients $a_k$ of $P(t)$.

For example, let $z_1$, $z_2$, and $z_3$ be the three roots of $P(t) \equiv x \cdot t^3 - (3+x) \cdot t^2 + (x^2 - 5x) \cdot t + 4$. Let $u_1 \equiv z_1 \cdot z_2^2 + z_3$ and let $u_2$ be a root of the polynomial $P_2(t) \equiv (7x^5 - 3x) \cdot t^2 + (6x - 1) \cdot t + 8$. Then there exists a minimal polynomial $P_3(t)$ over $\mathbb{Z}[x]$ whose roots are $u_1$ and $u_2$.

The polynomial $P_3(t)$ in Example 3.3 has a linear differential resolvent $\mathfrak{R}$ over $\mathbb{Z}[x]$. The terms of $\mathfrak{R}$ lie in the differential field generated by both the coefficients of $P(t)$ and $P_2(t)$, which happens to be the same field $\mathbb{Z}[x]$ for both $P(t)$ and $P_2(t)$. Hence, so $\mathfrak{R}$ is a joint linear differential resolvent *of* the polynomial $P(t)$ *for* the polynomial $y = z_1 \cdot z_2^2 + z_3$ and *of* the polynomial $P_2(t)$ *for* the polynomial $y = u_2$. Of course, it is much easier to say simply that $\mathfrak{R}$ is a resolvent of the polynomial $P_3(t)$.

## 4. Powersum formula generalized to resolvents of multiple polynomials for multivariable pseudopolynomials.

For the $N \cdot M$ (possibly indeterminate) constants $\{\alpha_{i,j}\}_{i \in [N], j \in [M]}$, the coefficient functions of the differential resolvent of the $N$ polynomials $P_i$ for the pseudopolynomial $y = \sum_{j=1}^{j=M} a_j \cdot \prod_{i=1}^{N} u_i^{\alpha_{i,j}}$ is of the form $\sum_{k=1}^{N} \sum_{s=1}^{h_k} r_{k,s} \cdot \prod_{i=1}^{L} \prod_{j=1}^{M} \alpha_{i,j}^{l_{i,j,k,s}} \cdot D^{m_k}$. As mentioned in Section 2 following Definition 2.9 and again in the final paragraph of the proof of Theorem 3.1, it

is understood by the subscript $k$ on the superscript $m_k$ and the subscript $s$ on $l_{i,j,k,s}$ that

$$\prod_{k=1}^{N}\prod_{s=1}^{h_k} r_{k,s} \neq 0.$$

**Definition 4.1.** Let $\mathbb{F}(\alpha)$ denote the differential field obtained by adjoining all the constants $\{\alpha_{i,j}\}_{i\in[N], j\in[M]}$ to the base differential field $\mathbb{F}$. Again, we assume $\mathbb{F}$ is the smallest differential field generated by the coefficients of the polynomials $\{P_i(t)\}_{i=1}^{N}$.

**Theorem 4.2.** If $\mathfrak{R} \equiv \sum_{k=1}^{N}\sum_{s=1}^{h_k} r_{k,s} \cdot \prod_{i=1}^{L}\prod_{j=1}^{M} \alpha_{i,j}^{l_{i,j,k,s}} \cdot D^{m_k}$ is a joint differential resolvent of the polynomials $\{P_i(t)\}_{i=1}^{L}$ over the differential field $\mathbb{F}(\alpha)$ for the pseudopolynomial

$$y \equiv \sum_{j=1}^{j=M} a_j \cdot \prod_{i=1}^{L} u_i^{\alpha_{i,j}}, \text{ with } \Psi \equiv \sum_{k=1}^{N} h_k \text{ nonzero coefficient functions } r_{k,s}, \text{ then we may}$$

compute the individual coefficient functions $r_{k,s}$ of $\mathfrak{R}$ by the *powersum formula*

$$r_{k,s} = \det\left[\left(\prod_{i=1}^{L}\prod_{j=1}^{M} \overline{\alpha}_{i,j}^{l_{i,j,k',s'}}\right) \cdot \left(\sum_{\sigma} D^{m_k}(\sigma\overline{y})\right)\right]_{\substack{(k',s')\neq(k,s) \\ (k',s')\times\overline{\alpha}}} \quad \text{(up to sign) } \textit{if} \text{ this formula yields a}$$

nonzero value. Notation will be explained in the proof.

**Proof.** Let $\overline{\alpha}$ represent one particular vector $\{\overline{\alpha}_{i,j}\}_{i\in[N], j\in[M]}$ of integers. Specialize $\alpha \equiv \{\alpha_{i,j}\}_{i\in[N], j\in[M]}$ to $\overline{\alpha}$, which we denote by $\alpha \to \overline{\alpha}$. When $\alpha \to \overline{\alpha}$, $y$ gets specialized to $\overline{y} \equiv \sum_{j=1}^{j=M} a_j \cdot \prod_{i=1}^{L} u_i^{\overline{\alpha}_{i,j}}$. If we specialize the equation

$$0 = \mathfrak{R}(\sigma y) = \sum_{k=1}^{N}\sum_{s=1}^{h_k} r_{k,s} \cdot \prod_{i=1}^{L}\prod_{j=1}^{M} \alpha_{i,j}^{l_{i,j,k,s}} \cdot D^{m_k}(\sigma y) \text{ by } \alpha \to \overline{\alpha}, \text{ we get}$$

$$0 = \mathfrak{R}(\sigma\overline{y}) = \sum_{k=1}^{N}\sum_{s=1}^{h_k} r_{k,s} \cdot \prod_{i=1}^{L}\prod_{j=1}^{M} \overline{\alpha}_{i,j}^{l_{i,j,k,s}} \cdot D^{m_k}(\sigma\overline{y}) \quad (4.1).$$

We sum (4.1) over all combined permutations (Definition 2.10) of the $\{P_i(t)\}_{i=1}^{L}$. We get

$$0 = \mathfrak{R}(\sigma\overline{y}) = \sum_{k=1}^{N}\sum_{s=1}^{h_k} r_{k,s} \cdot \prod_{i=1}^{L}\prod_{j=1}^{M} \overline{\alpha}_{i,j}^{l_{i,j,k,s}} \cdot \left(\sum_{\sigma} D^{m_k}(\sigma\overline{y})\right) \quad (4.2).$$

The sum $\sum_{\sigma} D^{m_k}(\sigma\overline{y})$ appearing in (4.2) is a multisum. Hence, one particular vector of positive integers $\overline{\alpha}$ yields one particular value for

$$\sum_{\sigma} D^{m_k}(\sigma\overline{y}) = \sum_{\sigma} D^{m_k} \sigma\left(\sum_{j=1}^{j=M} a_j \cdot \prod_{i=1}^{L} u_i^{\overline{\alpha}_{i,j}}\right) = \sum_{j=1}^{M} D^{m_k}\left(a_j \cdot \left(\sum_{\sigma}\prod_{i=1}^{L}(\sigma u_i)^{\overline{\alpha}_{i,j}}\right)\right) \quad (4.3).$$

The powersum formula got its name *powersum* because originally the author summed over *all* permutations of the roots of the polynomial $\{P_i(t)\}_{i=1}^{L=1} = \{P\}$, which yields the powersums of the roots of $P$. Although Definition 2.10 for a combined permutation $\sigma$ of the polynomials $\{P_i(t)\}_{i=1}^{L}$ was given in terms of *all* permutations of each of the roots of each $P_i(t)$, for Theorem 4.2 it is sufficient that $\sigma$ be given only in terms of the *Galois* permutations of the roots of each $P_i(t)$. So, we may expand the inner summation in (4.3) to show more clearly how the coefficients of the polynomials $\{P_i(t)\}_{i=1}^{L}$ arise. We have $\sum_{\sigma}\prod_{i=1}^{N}(\sigma u_i)^{\bar{\alpha}_{i,j}} = \sum_{\sigma_1}\sum_{\sigma_2}\cdots\sum_{\sigma_N}(\sigma_1 u_1)^{\bar{\alpha}_{1,j}} \cdot (\sigma_2 u_2)^{\bar{\alpha}_{2,j}} \cdots (\sigma_N u_N)^{\bar{\alpha}_{N,j}}$ where we sum over each member $\sigma_1$ of the Galois group $G_1$ of $P_1$, and sum over each member $\sigma_2$ of the Galois group $G_2$ of $P_2$, etc. So

$$\sum_{\sigma}\prod_{i=1}^{L}(\sigma u_i)^{\bar{\alpha}_{i,j}} = \left(\sum_{\sigma_1}(\sigma_1 u_1)^{\bar{\alpha}_{1,j}}\right)\cdot\left(\sum_{\sigma_2}(\sigma_2 u_2)^{\bar{\alpha}_{2,j}}\right)\cdots\left(\sum_{\sigma_N}(\sigma_N u_N)^{\bar{\alpha}_{N,j}}\right).$$ For each $i \in [L]$, $\left(\sum_{\sigma_i \in G_i}(\sigma_i u_i)^{\bar{\alpha}_{i,j}}\right)$ lies in the (non-differential) field generated by the coefficients of $P_i$.

Hence, for each $i \in [N]$, $\left(\sum_{\sigma_i \in G_i}(\sigma_i u_i)^{\bar{\alpha}_{i,j}}\right)$ is a multivariable polynomial over $\mathbb{Q}$ of the $\bar{\alpha}_{i,j}$ *powersums* of the roots of $P_i$. Hence, $\sum_{\sigma}D^{m_k}(\sigma\bar{y})$ lies in the smallest differential field differentially generated by the coefficients of all the $\{P_i(t)\}_{i=1}^{N}$ in (4.3) and (4.2).

Let $M_{k,s}$ denote the $(\Psi-1) \times (\Psi-1)$ matrix

$$\left[\left(\prod_{i=1}^{L}\prod_{j=1}^{M}\bar{\alpha}_{i,j}^{l_{i,j,k',s'}}\right)\cdot\left(\sum_{\sigma}D^{m_k}(\sigma\bar{y})\right)\right]_{\substack{(k',s')\neq(k,s)\\(k',s')\times\bar{\alpha}}}$$ whose rows are indexed by $(k,s)$, the index of the coefficient-function $r_{k,s}$, and whose columns are indexed by the particular vector of integers $\bar{\alpha}$. Then the system of equations (4.1), with each equation indexed by $\bar{\alpha}$, is a system of linear equations in the $r_{k,s}$. The solution is provided by Kramer's rule. Hence, $r_{k,s} = \det M_{k,s}$ (up to the appropriate sign). Since

$$\det\left[\left(\prod_{i=1}^{L}\prod_{j=1}^{M}\bar{\alpha}_{i,j}^{l_{i,j,k',s'}}\right)\cdot\left(\sum_{\sigma}D^{m_k}(\sigma\bar{y})\right)\right]_{\substack{(k',s')\neq(k,s)\\(k',s')\times\bar{\alpha}}}$$ lies in the smallest differential field generated by the $\{P_i(t)\}_{i=1}^{L}$, either all $r_{k,s} = 0$ or $\mathfrak{R}$ is a joint resolvent. **Q.E.D.**

**Complexity estimates of the powersum formula for joint resolvents.**

As example 1.4 showed, the number of coefficient-functions $r_{k,s}$ of a joint resolvent of two quadratic polynomials can be on the order of 2646. The minimal number of non-zero coefficient-functions for a joint resolvent may be much lower. This requires proof. The powersum formula for a resolvent with 2646 coefficient-functions would require the computation of the determinant of a 2645 x 2645 matrix with polynomial differential entries. Such a determinant would require on the order of 2645! sums and differences. Stirling's formula says this is on the order of $10^{7905}$ steps. This does not include the 2645 multiplications which must occur in each summand.

## 5. Difficulties in generalizing resolvents to arbitrary compositions of sums and powers

It is worth asking the question: does there exist a linear differential resolvent $\mathfrak{R}$ of the pseudo pseudopolynomial $y = (z^\alpha + z^\beta)^\lambda$ for constant indeterminates $\alpha$, $\beta$, and $\lambda$? To be consistent with earlier definitions of a linear differential resolvent of a polynomial, $P(t)$, the terms of $\mathfrak{R}$ must lie in the differential ring $\mathbb{F}[\alpha, \beta, \lambda]$. The author believes that this conjecture is false. We assume the coefficients of $t$ in $P(t)$ generate a differential field $\mathbb{F}$ over $\mathbb{Q}$.

We know that a resolvent exists for the pseudopolynomial $u = z^\alpha + z^\beta$ over the differential ring $\mathbb{F}[\alpha, \beta]$. But does that mean a resolvent for $u^\lambda$ exists over $\mathbb{F}[\alpha, \beta, \lambda]$?

Hence we may formulate this question slightly more generally. Suppose $u$ has a linear differential resolvent over some differential field $\mathbb{F}$. Does $u^\lambda$ have a linear differential resolvent over $\mathbb{F}[\lambda]$? Not necessarily. Suppose that a given differential equation $\sum_{k=0}^{n} h_k \cdot D^k u = 0$ of order $n$ has some fundamental system of $n$ solutions $\{u_i\}_{i=1}^{n}$. One can then *define* $e_k$ to be the $k$-th elementary symmetric function in the $\{u_i\}_{i=1}^{n}$. Then the polynomial $P(t) \equiv \sum_{k=0}^{n} e_k \cdot (-1)^{n-k} t^k$ depends upon the particular linear combination of the $\{u_i\}_{i=1}^{n}$ over constants. Thus, $P(t)$ depends upon the arbitrary constants chosen. The powersum formula will yield a resolvent $\sum_{k=0}^{n} F_k \cdot D^k u = 0$ with each $F_k \in \mathbb{Z}\{e_1, \ldots, e_n\}$. The differential field generated by the $h_k$ is *not* the same as the differential field generated by the $e_k$. Except for some common factor of the $h_k$, the differential field generated by the $h_k$ is a *proper* differential subfield of the differential field generated by the $e_k$. Hence, there exists some common factor $\Theta \neq 0$ such that $F_k = \Theta \cdot h_k$ and $\Theta \in \mathbb{Z}\{e_1, \ldots, e_n\}$.

In fact, the author is doing research using these ideas to determine the transcendentality of numbers of the form $u(1)$ when $\sum_{k=0}^{n} h_k \cdot D^k u = 0$ is given, $h_k \in \mathbb{Z}[x]$, and $D^k u(0)$ for each $k \in [n-1]_0$ is a specified integer. The simplest such example would be the differential equation $Du - u = 0$ with $u(0) = 1$ and the transcendentality of $u(1) = e$. The aim of this line of research is to generalize the proofs of the transcendentality of $e$ and $\pi$. The author was first made aware of the short 3-page proofs of the transcendentality of $e$ and $\pi$ at the ECCAD conference in [11].

# 6. Demonstrating the powersum formula Theorem 4.2 for a joint resolvent

The theorems in [1] come in very handy for setting bounds on the total degrees of the independent variables in resolvents.

**Example 6.1.** One of the simplest *computable* examples using Theorem 4.2 is the joint linear differential resolvent over $\mathbb{Z}[x]$ for the pseudopolynomial $v^{\sqrt{7}} + z^{\pi}$ for the linear polynomials $P_z(t) \equiv t - x$, $P_v(t) \equiv t - (x+1)$, that is, for $x^{\sqrt{7}} + (x+1)^{\pi}$. One can easily compute this resolvent entirely by hand. To make use of Theorem 4.2, however, let us find the form of the resolvent for the pseudopolynomial $y = v^{\alpha} + z^{\beta}$. We need differentiate $y$ only twice, $Dy = v^{\alpha} \cdot \dfrac{\alpha}{x} + z^{\beta} \cdot \dfrac{\beta}{x+1}$ and

$D^2 y = v^{\alpha} \cdot \dfrac{\alpha \cdot (\alpha - 1)}{x^2} + z^{\beta} \cdot \dfrac{\beta \cdot (\beta - 1)}{(x+1)^2}$, and eliminate $v^{\alpha}$ and $z^{\beta}$. Hence the joint resolvent for $y$ of $P_z$ and $P_v$ equals

$$\det \begin{bmatrix} y & 1 & 1 \\ Dy & \alpha/x & \beta/(x+1) \\ D^2 y & \alpha \cdot (\alpha - 1)/x^2 & \beta \cdot (\beta - 1)/(x+1)^2 \end{bmatrix} = 0. \tag{6.1}$$

Clearing denominators in (6.1) and multiplying by -1 yields

$$(\alpha \cdot (x+1)^2 \cdot x - \beta \cdot (x+1) \cdot x^2) \cdot D^2 y$$
$$+ (\alpha \cdot (x+1)^2 - \alpha^2 \cdot (x+1)^2 - \beta \cdot x^2 + \beta^2 \cdot x^2) \cdot Dy \tag{6.2}$$
$$+ (\alpha^2 \cdot \beta \cdot (x+1) - \alpha \cdot \beta^2 \cdot x - \alpha \cdot \beta) \cdot y = 0$$

So the joint resolvent has the form.

$$(r_{2,1} \cdot \alpha + r_{2,2} \cdot \beta) \cdot D^2 y$$
$$+ (r_{1,1} \cdot \alpha + r_{1,2} \cdot \alpha^2 + r_{1,3} \cdot \beta + r_{1,4} \cdot \beta^2) \cdot Dy \tag{6.3}$$
$$+ (r_{0,1} \cdot \alpha^2 \cdot \beta + r_{0,2} \cdot \alpha \cdot \beta^2 + r_{0,3} \cdot \alpha \cdot \beta) \cdot y = 0$$

Since $P_z$ and $P_v$ each have only 1 root, the sum of $y = v^{\alpha} + z^{\beta}$ over each of the roots of $P_z$ and then over each of the roots of $P_v$ means

$$\sum_{\substack{z \ni P_z(z)=0 \\ v \ni P_v(v)=0}} y = \sum_{\substack{z \ni P_z(z)=0 \\ v \ni P_v(v)=0}} (v^{\alpha} + z^{\beta}) = \sum_{v \ni P_v(v)=0} (v^{\alpha} + (x+1)^{\beta}) = x^{\alpha} + (x+1)^{\beta}.$$

Observe that in general

$$\sum_{\substack{z \ni P_z(z)=0 \\ v \ni P_v(v)=0}} y = \sum_{\substack{z \ni P_z(z)=0 \\ v \ni P_v(v)=0}} (v^{\alpha} + z^{\beta}) = \sum_{z \ni P_z(z)=0} (\tilde{p}_{\alpha} + \deg_v P_v \cdot z^{\beta}) = (\deg_z P_z) \cdot \tilde{p}_{\alpha} + (\deg_v P_v) \cdot \hat{p}_{\beta}$$

So Theorem 4.2 gives the formula for $r_{2,1}$, $r_{2,2}$, $r_{1,1}$, $r_{1,2}$, $r_{1,3}$, $r_{1,4}$, $r_{0,1}$, $r_{0,2}$, $r_{0,3}$ as the determinant of the 8 by 8 matrix formed from the following entries by setting $\alpha$ and $\beta$ to any set of 8 distinct pairs of integer values.

$$\alpha \cdot D^2(x^\alpha + (x+1)^\beta)$$
$$\beta \cdot D^2(x^\alpha + (x+1)^\beta)$$
$$\alpha \cdot D(x^\alpha + (x+1)^\beta)$$
$$\alpha^2 \cdot D(x^\alpha + (x+1)^\beta)$$
$$\beta \cdot D(x^\alpha + (x+1)^\beta) \tag{6.4}$$
$$\beta^2 \cdot D(x^\alpha + (x+1)^\beta)$$
$$\alpha^2 \cdot \beta \cdot (x^\alpha + (x+1)^\beta)$$
$$\alpha \cdot \beta^2 \cdot (x^\alpha + (x+1)^\beta)$$
$$\alpha \cdot \beta \cdot (x^\alpha + (x+1)^\beta)$$

The formula for $r_{2,1}$ omits the $\alpha \cdot D^2(x^\alpha + (x+1)^\beta)$ row.

The formula for $r_{2,2}$ omits the $\beta \cdot D^2(x^\alpha + (x+1)^\beta)$ row, and so forth.

In Maple, $(\alpha, \beta)$ was set to (1,1), (1,2), (1,3), (1,4), (2,1), (2,2), (2,3), (2,4). The Maple commands, written below only in ASCII (no Mathtype), were

```
with(LinearAlgebra):
Mrow:=[alpha*diff(diff(x^alpha + (x+1)^beta,x),x),
beta*diff(diff(x^alpha + (x+1)^beta,x),x),
alpha*diff(x^alpha + (x+1)^beta,x), alpha^2* diff(x^alpha + (x+1)^beta,x),
beta*diff(x^alpha + (x+1)^beta,x), beta^2* diff(x^alpha + (x+1)^beta,x),
alpha^2*beta*(x^alpha + (x+1)^beta),alpha*beta^2*(x^alpha +
(x+1)^beta),alpha*beta*(x^alpha + (x+1)^beta)];

#Choose letter M to stand for "matrix"
M:=Matrix([eval(Mrow,{alpha=1,beta=1}),eval(Mrow,{alpha=1,beta=2}),
eval(Mrow,{alpha=1,beta=3}), eval(Mrow,{alpha=1,beta=4}),
eval(Mrow,{alpha=2,beta=1}), eval(Mrow,{alpha=2,beta=2}),
eval(Mrow,{alpha=2,beta=3}), eval(Mrow,{alpha=2,beta=4})]);

#All the t(m,s) will be a common multiple of our desired final answer r(m,s)
#t for "test". The negative sign in front of factor comes from Kramer's rule from which
theorem 4.2 was derived.
t21:= factor(Determinant(DeleteColumn(M,1)));
t22:= -factor(Determinant(DeleteColumn(M,2)));
t11:= factor(Determinant(DeleteColumn(M,3)));
t12:= -factor(Determinant(DeleteColumn(M,4)));
t13:= factor(Determinant(DeleteColumn(M,5)));
t14:= -factor(Determinant(DeleteColumn(M,6)));
t01:= factor(Determinant(DeleteColumn(M,7)));
t02:= -factor(Determinant(DeleteColumn(M,8)));
t03:= factor(Determinant(DeleteColumn(M,9)));
```

The Maple output was

$$t_{2,1} = x \cdot (x+1)^2 \cdot \chi$$
$$t_{2,2} = -x^2 \cdot (x+1) \cdot \chi$$
$$t_{1,1} = (x+1)^2 \cdot \chi$$
$$t_{1,2} = -(x+1)^2 \cdot \chi$$
$$t_{1,3} = -x^2 \cdot \chi \quad (6.5)$$
$$t_{1,4} = x^2 \cdot \chi$$
$$t_{0,1} = -x \cdot \chi$$
$$t_{0,2} = (x+1) \cdot \chi$$
$$t_{0,3} = -\chi$$

where $\chi = 128 \cdot x \cdot \rho$ where
$$\rho = (840 \cdot x^9 + 4092 \cdot x^8 + 8212 \cdot x^7 + 8858 \cdot x^6 + 5554 \cdot x^5$$
$$+1899 \cdot x^4 + 159 \cdot x^3 - 144 \cdot x^2 - 63 \cdot x - 9)$$

Comparison of (6.5) with (6.2) and (6.3) shows that $t_{m,s} = \chi \cdot r_{m,s}$.

The simultaneously fascinating and frustrating aspect of using the powersum formula is the factors, such as $\chi$ and $\rho$, which arise and contribute to massive intermediate "blowup" difficulties. The Maple commands
with(RootFinding): Analytic(rho, x, re=-2..2,im=-2..2);
determine the nine roots of $\rho$. Three are real. Six are complex.

$$\rho \approx 840 \cdot (x + 0.7491227520) \cdot (x + 0.3964865679) \cdot (x - 0.2982846017)$$
$$\cdot (x^2 + 2.612786130 \cdot x + 1.733734989) \cdot (x^2 + 0.9607495874 \cdot x + 0.3670286507) \quad (6.6)$$
$$\cdot (x^2 + 0.4505681358 \cdot x + 0.1900501608)$$

None of the roots in (6.6) has an apparent relation with (6.2).

A more factored form of (6.3) allows us to see actual and apparent singularities of the resolvent more clearly.

$$((\alpha - \beta) \cdot x + \alpha) \cdot x \cdot (x+1) \cdot D^2 y$$
$$+((\alpha - \alpha^2 - \beta + \beta^2) \cdot x^2 + (\alpha - \alpha^2) \cdot (2x+1)) \cdot Dy \quad (6.7)$$
$$+((\alpha - \beta) \cdot x + (\alpha - 1)) \cdot \alpha \cdot \beta \cdot y = 0$$

We complete the solution to Example 6.1 by specializing $\alpha \to \sqrt{7}$ and $\beta \to \pi$ in (6.7)

$$((\sqrt{7} - \pi) \cdot x + \sqrt{7}) \cdot x \cdot (x+1) \cdot D^2 y$$
$$+((\sqrt{7} - 7 - \pi + \pi^2) \cdot x^2 + (\sqrt{7} - 7) \cdot (2x+1)) \cdot Dy \quad (6.8)$$
$$+((\sqrt{7} - \pi) \cdot x + (\sqrt{7} - 1)) \cdot \sqrt{7} \cdot \pi \cdot y = 0$$

A handheld calculator is sufficient to evaluate (6.8) numerically. The result, precise to as many digits as the HP48SX calculator would allow in each coefficient of $x$, is

$$(2.64575131106 - 0.49584134253 \cdot x) \cdot x \cdot (x+1) \cdot D^2 y$$
$$+(2.37376305856 \cdot x^2 - 8.708497378 \cdot x - 4.35424868894) \cdot Dy \quad (6.9)$$
$$+(-4.12137020879 \cdot x + 13.679275693) \cdot y = 0$$

**Example 6.2.** Choose the same problem as Example 6.1. However, this time, rather than deliberately choosing the smallest possible set of pairs of integer values for $(\alpha, \beta)$, let $(\alpha, \beta)$ run over the randomly chosen set (1,5), (2,3), (4,7), (1,4), (3,2), (8,1), (5,5), (6,4) in Theorem 4.2. Computation with Maple shows that the common "extraneous" factor $\chi$ from the powersum formula is $\chi = 41600 \cdot (x+1)^2 \cdot \rho$ where $\rho$ equals a polynomial with integer coefficients of degree 30 in $x$.

It is safe to conjecture that one may substitute any set $S$ of positive integer values for the indeterminates $\alpha_{i,j}$ into the powersum formula Theorem 4.2 and the formula will yield an answer which is not identically zero, provided of course that $|S|+1$ equals the number of undetermined coefficients $r_{m,s}$ in the resolvent. However, this result remains unproven even in the single polynomial (non-joint resolvent) case, except for the most "generic" case of a polynomial whose coefficients are all differentially transcendental over constants. See [11] for the proof of this result in this case.

**Acknowledgements.** The author wishes to thank Robert Israel of the University of British Columbia and John Fredsted of Aarhus University in Denmark on the MaplePrimes webboard for their assistance in locating the correct Maple commands. Finally, special thanks to Professor Roger Grimshaw, Director, Centre for Nonlinear Mathematics and Applications, of Loughborough University, Loughborough, UK and his two anonymous referees for reviewing my paper, making corrections, and urging me to expand my exposition.